\documentclass[a4paper,10pt]{article}

\usepackage[english]{babel}
\usepackage[T1]{fontenc}
\usepackage[dvips]{graphicx}
\usepackage{amsmath,amsfonts,amssymb,amsthm,bbm,latexsym,graphicx}
\usepackage{natbib}
\usepackage{a4wide}

\newtheorem{theo}{Theorem}

\newtheorem{lem}{Lemma}
 \theoremstyle{definition}
 \theoremstyle{remark}
\newtheorem{rem}{Remark} 
\numberwithin{equation}{section}

 \def\Zset{\mathbb{Z}} \def\Nset{\mathbb{N}}

\def\epsilon{\varepsilon}

\def\phi{\varphi}

\begin{document}

\title{Aggregation of isotropic autoregressive fields}
\author{Fr\'ed\'eric Lavancier${}$\\
{\footnotesize ${}$ Universit\'e de Nantes, Laboratoire de Mathématiques Jean Leray}\\
{\footnotesize  UFR Sciences et Techniques}\\{\footnotesize 2 rue de la Houssinière - BP 92208 -}\\ {\footnotesize  F-44322 Nantes Cedex 3,
France}}
\date{\null}

\maketitle

\vskip 3 cm
\begin{abstract}
This note constitutes a corrigendum to the article of \cite{azomahou}, JSPI, 139:2581-2597. The aggregation of isotropic four nearest neighbors autoregressive models on the lattice $\Zset^2$, with random coefficient, is investigated. The spectral density of the resulting random field is studied in details for a large class of law of the AR coefficient. Depending on this law, the aggregated field may exhibit short memory or isotropic long memory.

\end{abstract}
\newpage

\section{Introduction}

Contemporaneous aggregation of random parameter autoregressive processes is a well-known scheme to produce long memory  in time series. This idea first appeared in \cite{robinson78} and was more deeply investigated in \cite{gran801}. A lot of papers have then been devoted to this subject, among them \cite{Goncalvez88}, \cite{oppenheimviano1999}, \cite{zaffaroni}. This construction of long memory is famous in econometrics, since it illustrated the fact that long memory may arise in macroeconomic time series, although micro-dynamics are short range dependent.

More recently, the concept of long memory has been extended to random fields on a lattice. In this setting,  the intensity of long memory may depend on the direction or may be isotropic. We refer to \cite{lavspringer} for a review. From a statistical point of view, the isotropic case is investigated in \cite{anhlunney95}, \cite{frias08} and \cite{wang09}. Some  anisotropic models are studied in  \cite{boissy05}, \cite{guo09} and \cite{beran09}. A statistical test to detect long memory in random fields is provided in \cite{lav08}.

As for time series, it seems natural to study the effect of aggregation in random fields, in particular the emergence of long memory. The aggregation in this case in not contemporaneous but is implemented at each point of the lattice. This scheme has been described in \cite{lavspringer}, where it is applied to produce some anisotropic models of long memory. In \cite{azomahou}, aggregation of isotropic four nearest neighbors autoregressive models on the lattice $\Zset^2$ is studied. These autoregressive models are very natural models on $\Zset^2$ and their introduction dates back to \cite{whittle54}.  Their aggregation is expected to produce isotropic long memory in some cases, as announced in \cite{azomahou}. Unfortunately, this article contains false results and  the emergence of isotropic long memory by aggregation is not clearly described.

The present note is a corrigendum to the article of \cite{azomahou}. In section \ref{aggreg}, the general scheme of aggregation for random fields on $\Zset^2$ is recalled. In section \ref{iso}, the aggregation of isotropic four nearest neighbors autoregressive models is considered. First, we recall the necessary and sufficient conditions for existence of a stationary solution to such autoregressive equations, which are known since \cite{whittle54} (see e.g.  \cite{Guyon}), but are improperly stated in  Lemma 1 in  \cite{azomahou}. Then, the existence of a stationary aggregated field in $L^2$ is studied.  Under some general conditions  on the law of the autoregressive parameter, this existence is ensured in Theorem \ref{th}.  The behavior of the spectral density of the induced field is also presented. The memory properties of the aggregated field are then deduced. These results are a corrected version of Propositions 2 and 3 in \cite{azomahou}.

\section{Aggregating autoregressive random fields with random coefficients}\label{aggreg}
The definition of an autoregressive field in $\Zset^2$ is  very similar to the definition of an autoregressive time series (see \cite{brockwell:davis:1991}) and appears in \cite{whittle54}. Basic results are presented in  \cite{rosenblatt85} and \cite{Guyon}.

Let $L_1$ (resp. $L_2$) be the lag operator acting on the first (resp. second) index of the field $(X_{i,j})_{(i,j)\in\Zset^2}$, i.e.
$L_1 X_{i,j}=X_{i-1,j}$ and $L_2 X_{i,j}=X_{i,j-1}$.

The field $(X_{i,j})_{(i,j)\in\Zset^2}$ is an autoregressive field if there exists a white noise $(\epsilon_{i,j})_{(i,j)\in\Zset^2}$  in $L^2$ and a  complex function $P(z_1,z_2)=\sum_{(k,l)\in S}a_{k,l}z_1^k z_2^l$, where $S$ is a finite subset of $\Zset^2$ and $(a_{k,l})_{(k,l)\in S}$ are real coefficients, such that for all $(i,j)\in\Zset^2$
\begin{equation}\label{defar}P(L_1,L_2)X_{i,j}=\epsilon_{i,j}.\end{equation}
It is well known that such an equation admits a stationary solution if and only if $P(e^{i\theta_1},e^{i\theta_2})\not=0$ for all $(\theta_1,\theta_2)$ (see  \cite{whittle54}, \cite{rosenblatt85} or \cite{Guyon}).  In this case, the field $(X_{i,j})_{(i,j)\in\Zset^2}$ may be written
\begin{equation}\label{mainfini}X_{i,j}=\sum_{(k,l)\in\Zset^2} c_{k,l}\epsilon_{i-k,j-l}\end{equation}
where ($c_{k,l}$) are the coefficients involved in  the Laurent expansion of $P^{-1}(z_1,z_2)$ in the neighborhood of $\{|z_1|=1,|z_2|=1\}$ and satisfy $\sum |c_{k,l}|<\infty$.
Note that the series in \eqref{mainfini} converges almost surely  and in $L^2$ (see \cite{brockwell:davis:1991}, Proposition 3.1.1, for a proof in the time series setting which can be easily adapted to random fields).

Let us now turn to autoregressive fields with random coefficients. They are defined as in \eqref{defar} where the coefficients $(a_{k,l})_{(k,l)\in S}$ defining $P$ are random. Assuming that for almost every $(a_{k,l})$,  $P(e^{i\theta_1},e^{i\theta_2})\not=0$ for all $(\theta_1,\theta_2)$, we deduce easily that the representation \eqref{mainfini} (where the coefficients $(c_{k,l})$ are now random) holds almost surely. See \cite{oppenheimviano1999} for a justification in the time series setting which remains valid for random fields. However, for  \eqref{mainfini} to be true in $L^2$, it is necessary and sufficient that 
\begin{equation}\label{ass1}
E\int_{[-\pi,\pi]^2}|P(e^{i\lambda_1},e^{i\lambda_2})|^{-2}d\lambda_1d\lambda_2<\infty,
\end{equation}
where the expectation is with respect to the law of the coefficients  $(a_{k,l})_{(k,l)\in S}$ (see the same references). If \eqref{ass1} holds,  the spectral density of $(X_{i,j})_{(i,j)\in\Zset^2}$ is 
\begin{equation}\label{spectral}f_X(\lambda_1,\lambda_2)=\frac{\sigma_\epsilon^2}{4\pi^2} E |P(e^{i\lambda_1},e^{i\lambda_2})|^{-2},\end{equation}
where $\sigma_\epsilon^2$ denotes the variance of the white noise $(\epsilon_{i,j})_{(i,j)\in\Zset^2}$.

Finally, aggregation of random fields is implemented according to the classical scheme for time series. Consider $N$ independent copies of autoregressive fields with random coefficients as above, i.e. for $n=1,\dots,N$, for all  $(i,j)\in\Zset^2$,
\begin{equation}P^{(n)}(L_1,L_2)X^{(n)}_{i,j}=\epsilon^{(n)}_{i,j},\end{equation}
where $(\epsilon^{(n)})$ are  i.i.d. copies of  white noises (called idiosyncratic components) and the $P^{(n)}$'s come from i.i.d. copies of the random coefficients $(a_{k,l})_{(k,l)\in S}$. Aggregation of these $N$ fields defines the field $Z^{(N)}$ as, for all  $(i,j)\in\Zset^2$,
\begin{equation}\label{Z}
Z^{(N)}_{i,j}=\frac{1}{\sqrt N}\sum_{n=1}^N X^{(n)}_{i,j}.
\end{equation}
Under condition \eqref{ass1}, since the $(X^{(n)}_{i,j})_n$'s are i.i.d, the field $Z^{(N)}$ shares the same $L^2$ properties as the $(X^{(n)}_{i,j})_n$'s. In particular, its spectral density is \eqref{spectral}. When $N$ goes to infinity, from the central limit theorem, $Z^{(N)}$ tends to a Gaussian field $Z$. The aggregated field $Z$ is  thus a zero mean Gaussian field with spectral density \eqref{spectral}. Moreover, contrary to $Z^{(N)}$,  $Z$ is ergodic (see \cite{robinson78}).

\section{Isotropic long memory random fields as aggregated fields}\label{iso}
In this section we study the particular aggregation of isotropic four nearest neighbors autoregressive fields as in \cite{azomahou}. Unfortunately that article contains false statements. Our results correct them and provide further properties.

For each $n\in\Nset$, $X^{(n)}$ is the autoregressive field following the representation, for all $(i,j)\in\Zset^2$, 
\begin{equation}\label{isotropicar}X^{(n)}_{i,j}-\theta\left(X^{(n)}_{i-1,j}+X^{(n)}_{i+1,j}+X^{(n)}_{i,j-1}+X^{(n)}_{i,j+1}\right)=\epsilon^{(n)}_{i,j},\end{equation}
where $\epsilon^{(n)}$ is a white noise in $L^2$. This representation corresponds to \eqref{defar} where $P(L_1,L_2)=1-\theta (L_1+L_1^{-1}+L_2+L_2^{-1})$. When $\theta$ is not random, this model has been first considered in \cite{whittle54} and is studied for instance in \cite{Guyon}. 
 Contrary to the statement of Lemma 1 in \cite{azomahou}, we deduce from these references (or by easy computation) that the condition $|\theta|<\frac{1}{4}$ is necessary and sufficient so that \eqref{isotropicar} admits a stationary solution in $L^2$.

Let us now assume that $\theta$ is random and belongs almost surely to $[0,\frac{1}{4}[$ (see Remark \ref{thetanegative} for the case $\theta<0$). More precisely, let us suppose that the law of $\theta$ admits a density of probability  of the form  $\Phi(x)(\frac{1}{4}-x)^{\alpha}$, for $0\leq x\leq \frac{1}{4}$ , where $\alpha>-1$ and $\Phi$ is some function whose properties will be specified later.
According to \eqref{spectral}, if \eqref{ass1} holds, the spectral density of $X^{(n)}$ defined by \eqref{isotropicar}  is in this case
$$f_X(\lambda_1,\lambda_2)=\frac{\sigma_\epsilon^2}{4\pi^2} \int_0^{\frac{1}{4}}\frac{\Phi(x)(\frac{1}{4}-x)^{\alpha}}{(1-2x(\cos(\lambda_1)+\cos(\lambda_2)))^2}dx.$$

Assumption \eqref{ass1} is equivalent to the integrability of $f_X$ on $[-\pi,\pi]^2$. Conditions on $\alpha$ and $\Phi$ which  imply  \eqref{ass1} are deduced from the following lemma.

\begin{lem}\label{lem1}
Let $\alpha>-1$ and $\Phi$ be a bounded function on $[0,\frac{1}{4}]$, continuous at $\frac{1}{4}$ with $\Phi(\frac{1}{4})\not=0$. Consider, for any $\sigma_\epsilon^2>0$, the function $f$ defined on $[-\pi,\pi]^2$ by
\begin{equation}\label{f}
f:\ (\lambda_1,\lambda_2)\mapsto \frac{\sigma_\epsilon^2}{4\pi^2}\int_0^{\frac{1}{4}}\frac{\Phi(x)(\frac{1}{4}-x)^{\alpha}}{(1-2x(\cos(\lambda_1)+\cos(\lambda_2)))^2}dx.\end{equation}
Then,
\begin{itemize}
\item[(i)] when $\alpha>1$, $f$ is continuous at any point $(\lambda_1,\lambda_2)\in[-\pi,\pi]^2$,
\item[(ii)] when $-1<\alpha\leq1$, $f$ is continuous at any  $(\lambda_1,\lambda_2)\not=(0,0)$ in $[-\pi,\pi]^2$ and diverges to $+\infty$ when $(\lambda_1,\lambda_2)\to(0,0)$. Moreover, as $(\lambda_1,\lambda_2)\to(0,0)$,
\begin{itemize}
\item if $-1<\alpha<1$,  
\begin{equation}\label{feq1}f(\lambda_1,\lambda_2)\sim c_\alpha (\lambda_1^2+\lambda_2^2)^{\alpha-1}\end{equation}
where $c_\alpha=\frac{\sigma_\epsilon^2}{4\pi^2}16^{-\alpha}\Phi(\frac{1}{4})\int_0^{+\infty}\frac{u^{\alpha}}{(1+u)^2}du$.
\item if $\alpha=1$ and assuming $\Phi$  is a $\beta$-Hölder function ($\beta>0$),
\begin{equation}\label{feq2}f(\lambda_1,\lambda_2)\sim c_1 |\ln(\lambda_1^2+\lambda_2^2)|\end{equation}
 where $c_1=\frac{\sigma_\epsilon^2}{4\pi^2}16^{-1}\Phi(\frac{1}{4})$.\end{itemize}
\end{itemize}
\end{lem}

The proof of this lemma is postponed to the appendix. As a consequence of Lemma \ref{lem1}, \eqref{ass1} holds if $\alpha>0$ and $\Phi$ satisfies the assumptions in Lemma \ref{lem1}. In this case, the aggregated field $Z$, constructed as the limit of \eqref{Z} when $N\to\infty$, belongs to $L^2$ and its spectral density is \eqref{f}. 
The memory properties of $Z$ are deduced from the behavior of its spectral density. When it is continuous, then $Z$ exhibits short memory. When the spectral density is unbounded, then  the autocovariance function of $Z$ is not summable and $Z$ is said to exhibit long memory. These properties are summarized in the following theorem.  It provides the spectral properties improperly presented  in Propositions 2 and 3 of \cite{azomahou}.

\begin{theo}\label{th}
Let $X^{(n)}$, $n\in\Nset$, be a sequence of i.i.d isotropic autoregressive fields defined by \eqref{isotropicar}, where $\theta$ is random with probability density  $\Phi(x)(\frac{1}{4}-x)^{\alpha}$, $0\leq x\leq \frac{1}{4}$. Assume that $\alpha>0$ and that $\Phi$ is a bounded function on $[0,\frac{1}{4}]$, continuous at $\frac{1}{4}$ with $\Phi(\frac{1}{4})\not=0$. If $\alpha=1$, assume moreover that $\Phi$  is a $\beta$-Hölder function with $\beta>0$. Let $Z$ be the aggregated field obtained as the limit of \eqref{Z} when $N\to\infty$. Then $Z$ is a zero mean Gaussian random field in $L^2$ with spectral density \eqref{f}, where $\sigma_\epsilon^2$ denotes the variance of the white noises in \eqref{isotropicar}. In particular, \begin{itemize}\item if $\alpha>1$, $Z$ is short-range dependent, for its spectral density is continuous everywhere, \item if $0<\alpha<1$, $Z$ exhibits long memory and its spectral density behaves as \eqref{feq1} at 0, \item if $\alpha=1$,  $Z$ exhibits long memory and its spectral density behaves as \eqref{feq2} at 0.
\end{itemize}
\end{theo}
According to the results recalled in section \ref{aggreg}, the proof of this theorem is a direct consequence of Lemma \ref{lem1}.

\begin{rem} Following Definition 1 in \cite{lavspringer}, a random field over $\Zset^2$ exhibits isotropic long memory if its
 spectral density behaves at 0 as 
$|\lambda|^{\gamma}L\left(\frac{1}{|\lambda|}\right)b\left(\frac{\lambda}{|\lambda|}\right)$, $
-2<\gamma<0$,  where $\lambda=(\lambda_1,\lambda_2)$, $|.|$ denotes the euclidean norm,  $L$ is slowly varying at infinity and where $b$ is continuous on the unit sphere in $\mathbb{R}^2$. According to \eqref{feq1}, this is the case for $Z$ when $0<\alpha<1$ with $\gamma=2\alpha-2$. When $\alpha=1$, the behavior of the spectral density at $0$ is also a diverging function of $|\lambda|$. For this reason, one may also say that $Z$ exhibits isotropic long memory, even if it does not suit the latter definition. 
\end{rem}

\begin{rem}\label{thetanegative}
From a general point of view,  the law of $\theta$ may be supported on the whole interval $]-\frac{1}{4},\frac{1}{4}[$. From Theorem \ref{th}, the behavior of the spectral density at $0$ is driven by the behavior of the probability density of $\theta$ near $\frac{1}{4}$. It is easy to see that under similar conditions, the behavior
of the probability density of $\theta$ at $-\frac{1}{4}$ rules the behavior of the spectral density at  points $(\lambda_1,\lambda_2)$ such that $|\lambda_1|=|\lambda_2|=\pi$. In particular, some similar spectral singularities may arise at these points, leading to the so-called seasonal long memory phenomena (see \cite{oppenheimouldhayeviano2000} in the time-series setting).
\end{rem}

\appendix

\section{Proof of Lemma \ref{lem1}}
For any fixed $x<\frac{1}{4}$, the integrand in \eqref{f} is continous with respect to $(\lambda_1,\lambda_2)\in[-\pi,\pi]^2$. Moreover, since $\cos(\lambda_1)+\cos(\lambda_2)\leq 2$ , it is bounded by $\frac{\Phi(x)}{4}(\frac{1}{4}-x)^{\alpha-2}$, which is integrable when $\alpha>1$. So $(i)$ is proved thanks to the dominated convergence theorem.

Similarly, for any $(\lambda_1,\lambda_2)$ in a sufficiently small neighborhood of some fixed $(\tilde\lambda_1,\tilde\lambda_2)\not=0$, there exists $\epsilon>0$ such that $\cos(\lambda_1)+\cos(\lambda_2)\leq 2-2\epsilon$. Thus for any fixed $x<\frac{1}{4}$, for any such $(\lambda_1,\lambda_2)$, the integrand in  \eqref{f}  is bounded by $\epsilon^{-2}\Phi(x)(\frac{1}{4}-x)^{\alpha}$. This gives the continuity of $f$ at any point  $(\tilde\lambda_1,\tilde\lambda_2)\not=0$ when $\alpha>-1$.

On the other hand, changing $\frac{1}{4}-x$ into $z$ in  \eqref{f} leads to 
$$f(\lambda_1,\lambda_2)= \frac{\sigma_\epsilon^2}{4\pi^2}\frac{4}{(2-\cos(\lambda_1)-\cos(\lambda_2))^2}\int_0^{\frac{1}{4}}\frac{|z|^\alpha \Phi(\frac{1}{4}-z)}{\left(1+4z\frac{\cos(\lambda_1)+\cos(\lambda_2)}{2-\cos(\lambda_1)-\cos(\lambda_2)}\right)^2}dz.$$

Letting $A_\lambda=\frac{\cos(\lambda_1)+\cos(\lambda_2)}{2-\cos(\lambda_1)-\cos(\lambda_2)}$ and applying the change of variables $u=4zA_\lambda$, we obtain
\begin{equation}\label{ftot}
f(\lambda_1,\lambda_2)= \frac{\sigma_\epsilon^2}{4\pi^2}\frac{4^{-\alpha}}{(2-\cos(\lambda_1)-\cos(\lambda_2))^{1-\alpha}}\frac{|\cos(\lambda_1)+\cos(\lambda_2)|^{-\alpha}}{\cos(\lambda_1)+\cos(\lambda_2)} \int_0^{A_\lambda} \frac{|u|^{\alpha} \Phi(\frac{1}{4}(1-\frac{u}{A_\lambda}))}{(1+u)^2}du.\end{equation}

Note that $A_\lambda\to+\infty$ as $(\lambda_1,\lambda_2)\to(0,0)$. Since $\Phi$ is continous at $\frac{1}{4}$ and bounded on $[0,\frac{1}{4}]$, the dominated convergence theorem applies when $-1<\alpha<1$ and easy computations give \eqref{feq1}.

When $\alpha=1$, we get from \eqref{ftot},
\begin{equation}\label{f_1}f(\lambda_1,\lambda_2)= \frac{\sigma_\epsilon^2}{4\pi^2}\frac{4^{-1}}{(\cos(\lambda_1)+\cos(\lambda_2))^{2}}\left(I_1+I_2\right),\end{equation} where

$$I_1=\int_0^{A_\lambda} \frac{|u|}{(1+u)^2}\Phi\left(\frac{1}{4}\right)du \quad \text{and} \quad I_2=\int_0^{A_\lambda} \frac{|u|}{(1+u)^2} \left[\Phi\left(\frac{1}{4}\left(1-\frac{u}{A_\lambda}\right)\right)-\Phi\left(\frac{1}{4}\right)\right]du.$$

Assuming $\Phi$ is a $\beta$-Hölder function on $[0,\frac{1}{4}]$ when $\alpha=1$, we have
\begin{equation}\label{I2}|I_2|<\frac{c}{|4A_\lambda|^{\beta}}\int_0^{A_\lambda} \frac{|u|^{1+\beta}}{(1+u)^2} du,\end{equation}
where $c>0$, which shows that $I_2\to 0$ as $(\lambda_1,\lambda_2)\to(0,0)$. For $I_1$, note that $A_\lambda>0$ if $(\lambda_1,\lambda_2)$ is small enough, so direct computation leads  to 
\begin{equation}\label{I1}I_1=\Phi\left(\frac{1}{4}\right)\left[\ln(1+A_\lambda)+\frac{1}{1+A_\lambda}-1\right].\end{equation}
The equivalence \eqref{feq2} then follows easily from \eqref{f_1}, \eqref{I2} and \eqref{I1}.

\bibliographystyle{apalike}
\bibliography{bibaph}

\end{document}